\documentclass[psfig]{article}

\usepackage{amssymb, amsmath, amsthm}

\newtheorem{theo}{Theorem}[section]
\newtheorem{lem}{Lemma}[section]

\newcommand{\ov}{\overline}

\newcommand{\BB}[1]{\ensuremath{\mathbb{#1}}}

\newcommand{\R}{\ensuremath{\BB{R}}}
\newcommand{\Z}{\ensuremath{\BB{Z}}}

\newcommand{\be}{\begin{equation}}
\newcommand{\ee}{\end{equation}}
\newcommand{\bes}{\begin{equation*}}
\newcommand{\ees}{\end{equation*}}
\newcommand{\bi}{\begin{itemize}}
\newcommand{\ei}{\end{itemize}}
\newcommand{\bea}{\begin{eqnarray}}
\newcommand{\eea}{\end{eqnarray}}
\newcommand{\beas}{\begin{eqnarray*}}
\newcommand{\eeas}{\end{eqnarray*}}

\begin{document}
\title{Multiplicity and regularity of large periodic solutions with rational frequency for a class of semilinear monotone wave equations.}
\author{Jean-Marcel Fokam \\
{\small \tt fokam\,@aun.edu.ng}\\
School of Arts and Sciences, American University of Nigeria\\
Yola, 
}


\date{} 



\maketitle

\begin{abstract}
We prove the existence of infinitely many classical large periodic solutions for a class of semilinear wave equations with periodic boundary conditions:
\[
u_{tt}-u_{xx}=f(x,u),
\]
\[
 u(0,t)=u(\pi,t)\,\,\,, u_x(0,t)=u_x(\pi,t).
\]
Our argument relies on some new estimates for the linear problem with periodic boundary conditions by combining Littlewood-Paley techniques, the Hausdorff-Young theorem of harmonic analysis, and a variational formulation due to Rabinowitz \cite{R78},\cite{R84}. We also develop a new approach to the regularity of the distributional solutions by differentiating the equations and employing Gagliardo-Nirenberg estimates .
\end{abstract}

\newpage

\section{Introduction}

\noindent In this paper we construct infinitely many large classical time-periodic solutions
for the following semilinear wave equation 
\begin{equation}
    u_{tt}-u_{xx}-f(x,u)=0. \label{papier}
\end{equation}
\be
    u(0,t)=u(\pi,t),\,\ u_x(0,t)=u_x(\pi,t)
\label{bdry}
\ee
where $f$ is $C^{2,1}$ and has polynomial growth and depends on $x,u$. The existence of large periodic solutions with periodic boundary conditions is not well understood. When the frequency is irrational the method of Craig and Wayne in \cite{CW}, extended to higher dimension by Bourgain\cite{Bourgain95} and Berti and Bolle \cite{BertiBolle2010}proves the existence of small periodic solutions for typical potentials but the existence of classical periodic solutions for rational frequency is not known. Note that typical constant potentials in \cite{CW},\cite{Bourgain95}, \cite{BertiBolle2010} are satisfyied for
\be
   u_{tt}-\Delta u-m-f(x,u)=0
\ee
for \textit{typical m}.
On the other hand there exists a substantial amount of literature for semilinear wave equations with Dirichlet boundary conditions see for instance \cite{R84},\cite{R78}, \cite{BCN} for rational frequencies and the proofs of existence of \textit{classical solutions with $f$ having some spatial dependence rely on a fundamental solution discovered by Lovicarova} in \cite{Lovicarova}. For the case of irrational frequencies the problem of periodic solutions with Dirichlet boundary conditions for the resonant case $m=0$ was shown by Lidskii and Schulman \cite{LidskiiShulman}, Berti and Bolle \cite{BertiBolleDuke}, \cite{BertiBolle2003} and Yuan \cite{Yuan} and, for some non zero constant potential in P\"oschel \cite{Poschel}, and typical potential by Kuksin \cite{Kuksin}. De Simon and Torelli in \cite{DST} do not employ Lovicarova's formula but their $C^0$ estimate relies on $L^2$ a priori estimates on $f(x,u)$ which are not readily avaliable for distributional solutions of (\ref{papier}). The difficulty in proving regularity of distributional solution of (\ref{papier}) stems for the kernel of $\Box$ which is infinite dimensional. In absence of such a fundamental solution for the d'Alembertian under periodic boundary conditions problem we develop an approach based on tools from harmonic analysis such as Littlewood-Paley techniques, the Hausdorff-Young theorem and Gagliardo-Nirenberg estimates. The Hausdorff-Young theorem had been employed earlier by Willem in \cite{Willem85} to get a $L^\infty$  a priori estimate on solutions and by Coron to prove a Sobolev embedding in \cite{Coron}. In this paper the existence of \textit{classical solution for time periodic solutions with periodic boundary conditions of the semilinear wave equation (\ref{papier}) will be shown by proving the stronger $C^\gamma$ H\"older estimates than the $L^\infty$ in} \cite{Willem85}, and our approach gives an alternative proof of the existence of classical periodic solutions in the case of Dirichlet boundary conditions with semilinear term with some spatial dependence for $f(x,u)$ sufficiently smooth.
 
In section 1 we prove the linear estimates we need to prove the regularity of the solution. In section 2 we follow the scheme of \cite{R84} and \cite{R78} to construct a weak solutions and in section 3 we show the regularity of the solution by repeated differentiation of the equations and the linear estimates proved in section 1 and Gagliardo-Nirenberg inequalities .

Since our proof is of variational nature it is natural to ask if there is a notion of critical exponent or critical growth for this equation. An open question is then whether there are semilinear terms $f(x,u)$ of say exponential or super exponential type (as this paper deals with seminilear terms of polynomial type) for which there are large amplitude distributional solutions which are not classical ($f(x,u)$ being assumed to be smooth).

Our second contribution is a new Gagliardo-Nirenberg type inequality for the natural function space associated to the existence of periodic solution for the d'Alembertian. The estimate we prove allow us to recover an estimate in Coron \cite{Coron} by proving here a stronger estimate which also shows that the Sobolev embedding in \cite{Coron} is \textit{compact}. The proof follows by adapting the Fourier approach to the Sobolev embedding as in the paper by Chemin \cite{Chemin}.
We seek time-periodic solutions satisfying periodic boundary conditions so we seek functions $u\in\R$ with expansions of the form
\[
u(x,t)=\sum_{(j,k)\in\Z\times\Z}\widehat{u}(j,k)e^{i2jx}e^{ikt}
\]
and define the function space $E$:
\[
||u||^2_{E}=\sum_{2j\neq \pm k}\frac{|Q|}{4}|k^2-4j^2||\widehat{u}(j,k)|^2+\sum_{2j=\pm k}|4j^2||\widehat{u}(j,k)|^2+|\widehat{u}(0,0)|^2
\]
where $Q=[0,\pi]\times[0,2\pi]$ and define the functions spaces $E^+,E^-,N$ as follows:
\[
N=\{ u\in E, \,\, \widehat{u}(j,k)=0 \,\, {\rm{for}} \,\, 2|j|\neq |k|\}
\]
Note that in the case of periodic boundary conditions the structure of the kernel $N$ of $\Box$ is slightly different than in the case of Dirichlet boundary conditions. Here $v\in N$ we have
\be
v(x,t)=p_1(x+t)+p_2(x-t)
\ee
where the $p_1,p_2\in H^1(0,\pi)$ are arbitrary $\pi$-periodic functions.
\[
E^+=\{ u\in E, \,\, \widehat{u}(j,k)=0 \,\, {\rm{for}} \,\, |k|\leq 2|j|\}
\]
\[
E^-=\{ u\in E, \,\, \widehat{u}(j,k)=0 \,\, {\rm{for}} \,\, |k|\geq 2|j|\},
\]
$u=v+w$, $w=w^++w^-$ where $w\in E$, $w^+\in E^+$,$w^-\in E^-$ and $v\in N$  
and define the norm on $E\oplus N$
\[
||u||_{\beta,E}^2=||w^+||^2_E+||w^-||^2_E+\beta||v||^2_{H^1}.
\]
\be
I_\beta(u)=\int_Q[\frac{1}{2}(u_t^2-u_x^2-\beta (v^2+v_t^2)-F(x,u)]dxdt.
\ee
When $u$ is trigonometric polynomial, $I_\beta$ can also be represented in $E^m\oplus N^m$ as:
\[
I_\beta(u)=\frac{1}{2}(||w^+||_E^2-||w^-||_E^2-\beta(||v||^2_{L^2}+||v_t||^2_{L^2})-\int_Q F(x,u)dxdt .
\]
which can also be written as
where $\frac{\partial F(x,u)}{\partial u}=f(x,u)$, and first seek weak solution of the modified equation:
\be
\Box u=\beta v_{tt}-f(x,u)-\beta v
\label{mpapier}
\ee
and then send the parameter $\beta$ to zero.\\
{\bf Assumptions on $f(u)$:}\\
we assume that there are positive constants $c_0^1 \leq c_0^2,c^1_1,c_1^2$ such that
\be
c_0^1|u|^{s-1}u+c_1^1\leq f(x,u) \leq c_0^2|u|^{s-1}u+c_2^1
\label{conditions}
\ee
with $c_0^1>\frac{c_0^2}{s+1}$. These assumptions are satisfied by some nonlinearities of polynomial type. $f(x,u)$ must also be strongly monotone increasing:
\be
\frac{\partial f(x,u)}{\partial u}\geq\alpha>0
\label{conditions2}
\ee

\begin{theo}Under assumptions (\ref{conditions}),(\ref{conditions2}) and $f\in C^{2,1}$ the equation (\ref{papier}) admits infinitely many classical solutions.
\end{theo}
\section{Estimates}
Define $l^q=\{\hat{u}(j,k) s.t. \sum_{(j,k)\in\Z\times\Z}|\hat{u}(j,k)|^q <+ \infty\}$. 
\begin{theo}
\label{Riesz}The function $u=\sum_{2j\neq\pm k}\hat{u}(j,k)e^{2ijx+ikt}\in C^{\gamma}$ where $\gamma<1-\frac{1}{p}$ 
if
\be
\hat{u}(j,k)=\frac{\hat{f}(j,k)}{4j^2-k^2}
\ee
for $2j\neq \pm k$, 
$\hat{f}\in l^q$ and $\frac{1}{p}+\frac{1}{q}=1$.
\end{theo}
Proof:\\
Let $B_m$ the set
\[
B_m=\{ (j,k)\in\Z\times\Z \,\, 2|j|+|k|\leq 2.2^m \}
\]
and $\Delta_m$
\[
\Delta_m=B_m\setminus B_{m-1}
\]
so we have in $\Delta_m$
\[
2^m\leq 2|j|+|k|\leq 2.2^m
\]
and the $C^\gamma$ norm will be estimated by
\[
\sup_m 2^{\gamma m}||\Delta_m||_{C^0}.
\]
see \cite{Schlag} or \cite{Katznelson}. 
\begin{eqnarray}2^{\gamma m}||{\Delta}_m||_{C^0}
& = & 2^{\gamma m}||\sum_{(j,k)\in{\Delta}_m}\hat{u}(j,k)e^{i2jx}e^{ikt}||_{C^0}\nonumber\\
& = & ||\sum_{(j,k)\in\Delta_m}\frac{2^m\hat{f}(j,k)}{4j^2-k^2}e^{i2jx}e^{ikt}||_{C^0}\nonumber\\
& \leq & [\sum_{(j,k)\in\Delta_m}\frac{2^{\gamma mp}}{(|2j|^2-k^2)^p}]^{\frac{1}{p}}[\sum_{(j,k)\in \Z\times\Z}|\hat{f}(j,k)|^q]^{\frac{1}{q}}\nonumber\\
& \leq & [\sum_{(j,k)\in\Delta_m}\frac{(2|j|+|k|)^{\gamma p}}{(2|j|+|k|)^p(|2j|-|k|)^2)^p}]^{\frac{1}{p}}[\sum_{(j,k)\in \Z\times\Z}|\hat{f}(j,k)|^q]^{\frac{1}{q}}\nonumber\\
& \leq & [\sum_{(j,k)\in\Delta_m}\frac{1}{(2|j|+|k|)^{(1-\gamma) p}(|2j|-|k|))^p}]^{\frac{1}{p}}[\sum_{(j,k)\in \Z\times\Z}|\hat{f}(j,k)|^q]^{\frac{1}{q}}\nonumber\\
& \leq & c||\hat{f}||_{l^q}\leq c||f||_{L^p}
\end{eqnarray}
as long as $\gamma<1-\frac{1}{p}$ and the last inequality follows from the Hausdorff-Young theorem. \\
Remark: The argument here provides an alternate proof of the H\"older continuity of weak solutions of $\Box w=f$ where $f\in L^p\cap N^\perp$ where $N^\perp$ denotes the weak orthogonal of the kernel of $\Box$ with Dirichlet boundary conditions, proved in \cite{BCN} via Lovicarova's fundamental solution, for $1<p\leq 2$.\\
In the case that $p=2$ we have $u\in C^{0,\gamma}$ or similarly $f\in H^\alpha$ implies $u\in C^{\alpha+\frac{1}{2}}$.\\ Define
\be
u_{h_1,h_2}(x,t)=u(x+h_1,t+h_2)
\ee
and 
\be
\Delta_m^{++}=\{(j,k)\in \Z\times\Z \,(j,k)\in\Delta_m \,\, j\geq 0, k\geq 0\}
\ee
\be
\Delta_m^{+-}=\{(j,k)\in \Z\times\Z \, (j,k)\in\Delta_m \,\, j\geq 0, k<0\}
\ee
\be
\Delta_m^{-+}=\{ (j,k)\in \Z\times\Z \, (j,k)\in\Delta_m \,\, j<0, k\geq 0 \}
\ee
\be
\Delta_m^{--}=\{ (j,k)\in \Z\times\Z \, (j,k)\in\Delta_m \,\, j\leq 0, k< 0\}
\ee
and define $u^{--},u^{++},u^{-+},u^{+-}$ as 
\begin{eqnarray}\hat{u}^{++}(j,k) & = & \hat{u}(j,k) \,\ {\rm{if}} \,\ j\geq 0,k\geq 0\nonumber\\
                            & = & 0 \,\ {\rm{otherwise}}
\end{eqnarray}
\begin{eqnarray}\hat{u}^{+-}(j,k) & = & \hat{u}(j,k) \,\ {\rm{if}} \,\ j\geq 0,k< 0\nonumber\\
                            & = & 0 \,\ {\rm{otherwise}}
\end{eqnarray}
\begin{eqnarray}\hat{u}^{-+}(j,k) & = & \hat{u}(j,k) \,\ {\rm{if}} \,\ j< 0,k\geq 0\nonumber\\
                            & = & 0 \,\ {\rm{otherwise}}
\end{eqnarray}
\begin{eqnarray}\hat{u}^{--}(j,k) & = & \hat{u}(j,k) \,\ {\rm{if}} \,\ j< 0,k< 0\nonumber\\
                            & = & 0 \,\ {\rm{otherwise}}
\end{eqnarray}
\begin{theo} If $u^{++}\in C^{0,\gamma}$ then $u^{++}\in H^{\gamma^\prime}$ if $\gamma^\prime<\gamma$
\end{theo}
The analogue is also true for $u^{+-},u^{-+},u^{--}$.\\
Proof:\\
\begin{eqnarray}||u^{++}||_{H^{\gamma^\prime}}^2
& =    & \sum_m\sum_{(j,k)\in\Delta_m^{++}}(|2j|+|k|)^{2\gamma^\prime}|\hat{u}(j,k)|^2\nonumber\\
& \leq & \sum_m\sum_{(j,k)\in\Delta_m^{++}}2^{2(m+1)\gamma^\prime}|\hat{u}(j,k)|^2\nonumber\\
& \leq & \sum_m 2^{(m+1)\gamma^\prime}\sum_{(j,k)\in\Delta_m^{++}}|e^{i(2|j|+|k|)h}-1|^2|\hat{u}(j,k)|^2\nonumber
\end{eqnarray}
with $h=h(m)=\frac{2\pi}{3}2^{-m}$. \\
Remark: in the next line the sum over $\Delta_m^{++}$ is extended to the whole series but $h$ still depends on $m$. This is possible because of Parseval. \\Then
\begin{eqnarray}||u^{++}||_{H^{\gamma^\prime}}^2
& \leq & \sum_m2^{(m+1)\gamma^\prime}||u_{h(m),h(m)}-u||_{l^2}\\
\label{Uh}
& \leq & \sum_m2^{(m+1)\gamma^\prime}||u_{h(m),h(m)}-u||_{C^0}\nonumber\\
& \leq & \sum_m2^{(m+1)\gamma^\prime}||u||_{C^0,\gamma}^2|h(m)|^{2\gamma}\nonumber\\
& \leq & \sum_m2^{(m+1)\gamma^\prime}||u||_{C^0,\gamma}^2|\frac{2\pi}{3}2^{-m}|^{2\gamma}\nonumber\\
& \leq & c||u||_{C^0,\gamma}^2\sum_m2^{m(\gamma^\prime-\gamma)}\nonumber\\
& \leq & c(\gamma-\gamma^\prime)||u||^2_{C^{0,\gamma}}\nonumber
\end{eqnarray}
The estimates for $u^{--},u^{-+},u^{+-}$ follow similarly by re placing $u_{h,h}$ in the preceding argument (\ref{Uh}) by $u_{-h,-h}$,$u_{-h,h}$,$u_{h,-h}$.
We can conclude by noting:
\be
||u||^2_{H^\gamma}=||u^{++}||_{H^\gamma}^2+||u^{+-}||_{H^\gamma}^2+||u^{-+}||_{H^\gamma}^2+||u^{--}||_{H^\gamma}^2
\ee
\begin{lem}
\label{H1}
Let $f,w\in L^2(Q)$ such that 
\be
\widehat{f}(j,k)=0=\widehat{w}(j,k)\,\ for  \,\, 2j=\pm k
\ee
 and 
\be
(-k^2+4j^2)\widehat{w}(j,k)=\widehat{f}(j,k),\,\ for \,\, 2j\neq\pm k \,\, 
\ee
then $w\in H^1$
\end{lem}
Proof:
\begin{eqnarray}||w||_{H^1}^2
& = & \sum_{2j\neq \pm k}\frac{4j^2+k^2}{|4j^2-k^2|^2}|\hat{f}(j,k)|^2\nonumber\\
& = & \sum_{2j\neq \pm k}\frac{1}{2}\frac{(2j-k)^2+(2j+k)^2}{(2j-k)^2(2j+k)^2}|\hat{f}(j,k)|^2\nonumber\\
& \leq & \sum_{2j\neq \pm k}|\hat{f}(j,k)|^2\nonumber\\
& \leq & ||f||_{L^2}^2
\end{eqnarray}
Let $E^s$ be the closure of $\{e^{i2jx+ikt},\,\ 2j\neq \pm k\}$ under the norm
\[
||u||_{E^s}^2=\sum_{2j\neq\pm k}|\widehat{u}(j,k)|^2|k^2-4j^2|^s
\]
the we have the Sobolev estimate:
\begin{theo} $0<s<1$ the space $E^{s}$ is continuously embedded in $L^p$ where $p=\frac{2-s}{1-s}$.
\end{theo}
This theorem implies that the embedding in \cite{Coron} $E^1\subset L^p$ is compact, as $E^1\subset E^s$ is compact for $s<1$. We will show that it also implies a Gagliardo-Nirenberg inequality of the type:
\be
||u||_{L^p}\leq c(p)||u||_{L^2}^{1-s(p)}||u||_{E^1}^{s(p)}
\ee
Proof:\\
\be
f=f_{1,A}+f_{2,A}
\ee
where
\be
f_{1,A}=\sum_{2j\neq \pm k ,2|j|+|k|\leq A}\hat{f}(j,k)e^{i2jx}e^{ikt}
\ee
and
\be
f_{2,A}=\sum_{2j\neq \pm k ,2|j|+|k|>A}\hat{f}(j,k)e^{i2jx}e^{ikt}
\ee
\begin{eqnarray} |f_{1,A}| & \leq & \sum_{ 2j\neq \pm k ,2|j|+|k|\leq A}|\hat{f}(j,k)| \nonumber\\
                           & \leq & \sum_{2j\neq \pm k ,2|j|+|k|\leq A}|4j^2-k^2|^{-\frac{s}{2}}|4j^2-k^2|^{\frac{s}{2}}|\hat{f}(j,k)|
\end{eqnarray}
and applying Cauchy-Schwarz we have
\begin{eqnarray}|f_{1,A}|  & \leq & (\sum_{2j\neq \pm k ,2|j|+|k|\leq A}\frac{1}{|4j^2-k^2|^s})^{\frac{1}{2}}(\sum_{2j\neq \pm k ,2|j|+|k|\leq A}|4j^2-                                    k^2|^s|\hat{f}(j,k)|^2)^{\frac{1}{2}}\nonumber\\
                           & \leq & ||f||_{E^s}(\sum_{m,n\in\mathbb{N}\leq A}\frac{4}{m^sn^s})^{\frac{1}{2}}\nonumber\\
                           & \leq & c||f||_{E^s}(\int_1^A\frac{dm}{m^s}\int_1^{A}\frac{dn}{n^s})^{\frac{1}{2}}\nonumber\\
                           & \leq & c||f||_{E^s} A^{-s+1}.
\end{eqnarray}
Now we seek $A_\lambda$ such that 
\be
|f_{1,A}|\leq \frac{\lambda}{4}.
\label{chemin2}
\ee So we require the estimate
\be
c||f||_E A^{1-s}\leq \frac{\lambda}{4}
\label{chemin1}
\ee
this leads to the inequality
\[ A^{1-s}\leq \frac{\lambda}{4c||f||_{E^s}}\]
So let $A_\lambda$:
\[ A_\lambda= \frac{\lambda}{4c||f||_{E^s}}.  \]
Now 
\[ \int_{[0,\pi][0,2\pi]}|f(x,t)|^pdxdt=p\int_0^\infty y^{p-1}w(y)dy \]
where $w_f(y)=\{(x,t)\in [0,\pi][0,2\pi]: |f(x,t|>y\}$. Now $|f(x,t)|>\lambda$ implies $|f_{1,A}|>\frac{\lambda}{2}$ or
$|f_{2,A}|>\frac{\lambda}{2}$. Recalling (\ref{chemin2}) and the definition of $A_\lambda$ conclude that 
\be
|f_{2,A_\lambda}|>\frac{\lambda}{2}
\ee                       
and
\be
w_f(\lambda)\leq w_{f_{2,A_\lambda}}(\frac{\lambda}{2})
\ee
hence 
\begin{eqnarray}\int_{[0,\pi][0,2\pi]}|f(x,t)|^pdxdt & = & p\int_0^\infty \lambda^{p-1}w_f(\lambda)d\lambda \nonumber\\
                                                     & \leq & p\int_0^\infty \lambda^{p-1}w_{f_{2,A_\lambda}}(\frac{\lambda}{2})d\lambda\nonumber
\end{eqnarray}
since 
\be
w(\lambda)\leq \frac{1}{\lambda^2}\int_{|f|\geq \lambda}|f(x,t)|dxdt
\ee
\begin{eqnarray} \int_{[0,\pi][0,2\pi]}|f(x,t)|^pdxdt 
& \leq & \int_0^\infty\lambda^{p-3}\int_{|f_{2,A_\lambda}(x,t)|>\frac{\lambda}{2}}|f_{2,A_\lambda}(x,t)|^2dxdtd\lambda\nonumber\\
& \leq & \int_0^\infty\lambda^{p-3}\int_{[0,\pi][0,2\pi]}|f_{2,A_\lambda}(x,t)|^2dxdtd\lambda
\end{eqnarray}
Then we can invoke Parseval formula to deduce
\begin{eqnarray} \int_{[0,\pi][0,2\pi]}|f(x,t)|^pdxdt 
& \leq & \int_0^\infty\lambda^{p-3}\sum_{2j\neq\pm k}|\hat{f}_{2,A_\lambda}(j,k)|^2d\lambda\nonumber\\
& = & \int_0^\infty\lambda^{p-3}\sum_{2j\neq \pm k,2|j|+|k|>A_{\lambda}}|\hat{f}(j,k)|^2.
\label{f}
\end{eqnarray}
Now
\[ 2|j|+|k|\geq A_\lambda= \frac{\lambda}{4c||f||_{E^s}} \] 
implies
\[ \lambda \leq 4c||f||_{E^s}((2|j|+|k|))^{1-s}\]
We continue the estimate from (\ref{f}):
\begin{eqnarray} \int_{[0,\pi][0,2\pi]}|f(x,t)|^pdxdt 
& \leq & \sum_{2j\neq\pm k,2|j|+|k|\geq A_{\lambda}}\int_0^\infty\lambda^{p-3}|\hat{f}(j,k)|^2\nonumber\\
& \leq & \sum_{2j\neq\pm k,2|j|+|k|\geq A_{\lambda}}\int_0^{4c||f||_{E^s}(2|j|+|k|)^{1-s}}\lambda^{p-                                    3}|\hat{f}(j,k)|^2\nonumber\\
& \leq & \sum_{2j\neq\pm k,2|j|+|k|\geq A_{\lambda}}|\hat{f}(j,k)|^2\int_0^{4c||f||_{E^s}(2|j|+|k|)^{1-s}}\lambda^{p-                                    3}d\lambda\nonumber\\
& \leq & \sum_{2j\neq\pm k,2|j|+|k|\geq A_{\lambda}}|\hat{f}(j,k)|^2[\frac{\lambda^{p-2}}{p-2}]_0^{4c||f||_{E^s}(2|j|+|k|)^{1-s}}\nonumber\\
& \leq & \sum_{2j\neq\pm k,2|j|+|k|\geq A_{\lambda}}|\hat{f}(j,k)|^2\frac{1}{p-2}[4c||f||_{E^s}((2|j|+|k|))^{1-s}]^{p-2}\nonumber\\
\end{eqnarray}
now if $s=(1-s)(p-2)$ i.e. $s=\frac{p-2}{p-1}$ then
\be
\int_{[0,\pi][0,2\pi]}|f(x,t)|^pdxdt\leq \frac{(4c)^{p-2}}{p-2}||f||_{E^s}^p
\ee
and we have the following Gagliardo-Nirenberg inequality for $p>2$:
\be
||u||_{L^p}\leq c(p)||u||_{E^{s(p)}}\leq c(p)||u||_{L^2}^{1-s(p)}||u||_{E^1}^{s(p)}.
\ee
\section{Construction of the weak solution}
 For the Galerkin procedure we define the spaces:
\[
E^m=span\{\sin 2jx\cos kt, \sin 2jx\sin kt,\, \cos 2jx\cos kt,\, \cos 2jx\sin kt,  \,\,\, 2j+k\leq m, 2j\neq k \},
\]
\[
E^{-m}=span\{\sin 2jx\cos kt, \sin 2jx\sin kt,\, \cos 2jx\cos kt,\, \cos 2jx\sin kt, \,\,\, 2j+k\leq m \,\ 2j<k \},
\]
\[
E^{+l}=span\{\sin 2jx\cos kt, \sin 2jx\sin kt,\, \cos 2jx\cos kt,\, \cos 2jx\sin kt, \,\,\, 2j+k\leq l \,\ 2j>k \},
\]
\[
N^m=span\{ \sin 2jx\cos kt, \sin 2jx\sin kt,\, \cos 2jx\cos kt,\, \cos 2jx\sin kt, \,\,\ 2j\leq m \}
\]
which are employed in the minimax procedure. We denote by $P^m$ the projection of $E\oplus N$ into $E^m\oplus N^m$. The functional $I_\beta$ satisfies the Palais-Smale condition. The arguments follows as in \cite{R84}, we do not repeat the argument here.
\begin{lem}$\forall u\in E^{+l}$, there is a constant $C(l)$ independent of $\beta,m$ such that 
\be
I_\beta(u)\leq M(l)
\ee
\end{lem} 
Proof:\\
Let $u\in E^{+l}$
\begin{eqnarray}I_\beta(u)
& =    & \frac{1}{2}||w^+||_E^2-\frac{1}{2}||w^-||_E^2-\beta||v||_{H^1}^2-\int_QF(u)dxdt\nonumber\\
& \leq & \frac{1}{2}||w^+||_E^2-\frac{1}{2}||w^-||_E^2-\beta||v||_{H^1}^2-c(s)\int_Q\frac{|u|^{s+1}}{s+1}dxdt+d(f,s)\nonumber\\
& \leq & c(f,s)+\sup_{u\in E^{+n}}\frac{1}{2}||w^+||_E^2-c(s,Q)||u||_{L^2}^{s+1}
\end{eqnarray}
Now in $E^{+l}$
\be
||u||_E^2\leq l||u||_{L^2}^2
\ee
and on the other-hand
\be
\sup_{u\in E^{+l}}\frac{1}{2}||w^+||_E^2-c(s,Q)||u||_{L^2}^{s+1}>0
\ee
while as $||u||_E\rightarrow +\infty$ in $E^{+l}$ is dominated by $||u||_{L^2}^{s+1}$ as $s+1>2$ 
and is attained at say $\overline{u}$
hence we have 
\be
c(s,Q)||\overline{u}||_{L^2}^{s+1}\leq ||\overline{u}||_E^2\leq l||\overline{u}||_{L^2}^2
\ee
and we can conclude there is $M(l)$ depending on $l$ but independent of $\beta$ such that
\be
I_\beta(u)\leq M(l).
\ee
Also $E^{+l}$ is finite dimensional hence there is $R(l)$ such that for all $u\in E^{+l}\oplus E^{-m}\oplus N^m$ and $||u||_{E,\beta}\geq R(l)$ implies $I_{\beta}(u)\leq 0$.
\begin{theo} {\rm{Let}} $f$ be $C^1$, {\rm{for}} $l$ {\rm{large enough there is a distributional solution}} $u=v+w$ {\rm{of the modified problem}} (\ref{mpapier}) .
\end{theo}
Proof:\\
In this proof the constants may dependent on $\beta$ and $f$ but are independent of $m$.
The proof of this theorem here is slightly simpler from the one in \cite{R84} as we take advantage of the polynomial growth of the nonlinear term. We also and employ Galerkin approximation.\\
Let $u_l^m=w^m+v^m\in E^m\oplus N^m$ a distributional solution corresponding to the critical value $c_l$, and any $\phi\in E^m\oplus N^m$:
\be
I^\prime(u_l^m)\phi=0
\label{Galerkin1}\ee
now taking $\phi=v_{tt}^m\in N^m$ we have
\[
(\beta v_{tt}^m,v_{tt}^m)_{L^2}=(f(x,u_l^m),v_{tt}^m)_{L^2}+\beta(v_t^m,v_t^m)
\]
and by (\ref{conditions}) there are constant positive $c,d$ such that 
\[
\beta||v_{tt}^m||_{L^2}^2\leq c||u^s||_{L^2}||v_{tt}^m||_{L^2}+d||v_{tt}^m||_{L^2}
\]
\[
\beta||v_{tt}^m||_{L^2}\leq c||v_{tt}^m||_{L^2}
\]
hence
\[
||v_{tt}^m||_{L^2}\leq c(\beta)
\]
we now have
\[
w_{tt}^m-w_{xx}^m=\beta v_{tt}^m+P^mf(x,u^m_l)\in L^2
\]
hence 
$w^m\in H^1\cap C^{\gamma}$, $\gamma<\frac{1}{2}$
by theorem \ref{Riesz} and lemma \ref{H1}. This now implies $w^m\in H^2$, $w^m\rightarrow w$ as $m\rightarrow +\infty$ pointwise and $w\in H^1\cap C^{\gamma}$. Then if $\phi=v_{tttt}^m$ then
\[
(\beta v_{tt}^m,v_{tttt}^m)_{L^2}=(f(x,u_l^m),v_{tttt}^m)_{L^2}-\beta(v^m_t,v^m_{ttt})
\]
so there exists $c$ independent of $m$ 
\[
(\beta v_{ttt}^m,v_{ttt}^m)_{L^2}= (f_u(x,u_l^m)u_{lt}^m,v_{ttt}^m)_{L^2}-\beta(v^m_t,v_{ttt}^m)
\] 
and we deduce $||v_{ttt}^m||_{L^2}\leq c(\beta)$ hence $v^m_{ttt}\rightarrow v_{tt}\in C^0$ as $m\rightarrow+\infty$ hence $v$ is $C^2$ and $w$ is $C^\gamma$ by applying theorem \ref{Riesz} to (\ref{mpapier}) . We now have
\[
u_l^m\rightarrow u\in C^\gamma \,\, {\rm{as}}\,\ m\rightarrow +\infty
\]
and since (\ref{Galerkin1}) holds for any $\phi\in E^m\oplus N^m$ we can deduce
\be
I^\prime(u)\phi=0 \,\,\, \forall \phi\in E^m\oplus N^m,
\ee
now sending $m\rightarrow\infty$,  $u$ is a weak solution of (\ref{mpapier}).\\
Then we can define $g_\theta(u)=u(x,t+\theta)$. Define:
\be
{\cal G}=\{ g_\theta \,\ s.t. \,\, \theta\in [0,2\pi) \}
\ee
\be
V_l=N^m\oplus E^{-m}\oplus E^{+l}
\ee
\be
G_l=\{h\in C(V_l,E^m) \,\, {\rm{ such \ that}} \,\ h \,\, {\rm{satisfies}} \,\,\ \gamma_1-\gamma_4 \}
\ee
$Fix {\cal G}=\{u\in E \,\, s.t. \,\, g(u)=u \,\, \forall g \in {\cal G} \}=span \{ \cos{2jx},\,\, \sin 2jx,\,\  j\in\Z \}\subset E^-$. Define $P^{0m},P^{-m}$ the orthogonal projections from $E^m\oplus N^m$ onto respectively $N^m(=E^{0m}),E^{-m}$ and $P_l$ the orthogonal projection from $E^m\oplus N^m$ onto $V_l$.
\begin{displaymath}
\left\{ \begin{array}{l}
\gamma_1 \,\, h \,\ {\rm{is \ equivariant}} \\
\gamma_2 \,\, h(u)=u \,\ {\rm{if}} \,\ u\in Fix{\cal G}      \\
\gamma_3 \,\ {\rm{There exists }} r=r(h) \,\ h(u)=u \,\, {\rm{if}} \,\, u\in V_l\setminus B_{r(h)} \\
\gamma_4 \,\ u=w^++w^-+v\in V_l \,\ (P^{0m}+P^{-m})h(u)=\alpha(u)v+\alpha^-(u)w^-+\phi(u) \,\, {\rm{where}} \,\, \alpha,\ov{\alpha}\in C(V_l,[1,\ov{\alpha}]) \\
\end{array}\right.
\end{displaymath}
and $1<\ov{\alpha}$ depends on $h,\phi$ continuous.
Define
\be
c_l(\beta)=\inf_{h\in G_l}\sup_{u\in V_l}I_\beta(h(u))
\ee
and $c_l(\beta)\rightarrow +\infty$ as $l\rightarrow\infty$ independently of $m,\beta$ .
\begin{lem}
$c_l(\beta)\rightarrow +\infty$ as $l\rightarrow +\infty$ 
\end{lem}
Proof:\\
\begin{eqnarray}I_\beta(u)
& = & \frac{1}{2}||w^+||_E^2-\frac{1}{2}||w^-||_E^2-\beta||v||_{H^1}^2-\int_Q F(u)dxdt
\end{eqnarray}
and there exists by assumptions (\ref{conditions}) $c(s),d(s)>0$ such that
\begin{eqnarray}I_\beta(u)
& \geq & \frac{1}{2}||w^+||_E^2-\frac{1}{2}||w^-||_E^2-\beta||v||_{H^1}^2-               c(s)\int_Q|u|^{s+1}dxdt-d(s)\nonumber\\
\end{eqnarray}
and if $u\in \partial B_\rho\cap V_{l-1}^\perp$ we have
\begin{eqnarray}I_\beta(u)
& \geq & \frac{1}{2}||w^+||_E^2-c(s)\int_Q|u|^{s+1}dxdt-d(s)\nonumber\\
\end{eqnarray}
now by the Sobolev embedding there is $\theta(s)<1$ such that if $\widehat{u}(j,k)=0$ for $2j=\pm k$ we have
\be
||u||_{L^{s+1}}\leq ||u||_{E^{\theta(s)}}
\ee
hence
\begin{eqnarray}I_\beta(u)
& \geq & \frac{1}{2}||w^+||_E^2-c(s)(||u||_{L^2}^{1-\theta(s)}||u||_{E^1}^{\theta(s)})^{s+1}-                  d(s)\nonumber\\
& \geq & \frac{1}{2}\rho^2-\rho^{s+1}l^{(1-\theta(s))\frac{s+1}{s-1}}
\end{eqnarray}now if we choose a constant $C(s)$ large and $\rho=\frac{1}{C(s)}l^{(1-\theta(s))\frac{s+1}{s-1}}$
\begin{eqnarray}I_\beta(u)
& \geq & \frac{1}{4}\rho^2-d(s)\nonumber\\
\end{eqnarray}
now recalling by applying the corollary $2.4$ in \cite{FHR} to $P_{l-1}h\in C(\partial B_{\rho_l},V_{l-1})$
we have
\be
h(V_l)\cap\partial B_{\rho_l}\cap V_{l-1}^\perp\neq \o
\ee
then
\be
\sup_{V_l}I_\beta(h(u))\geq \inf_{u\in\partial B_{\rho_l}\cap V_{l-1}^\perp}I(u)\geq \frac{1}{4}\rho^2-d(s)\rightarrow +\infty
\ee
The $c_l(\beta)$ are critical values of $I_\beta$ on $E^m$. This is obtained by a standard argument see \cite{R84} propositions 2.33 and 2.37. 
\begin{lem}If $u$ is a critical point of $I_\beta$ in $E^m\oplus N^m$ then there are constants $c_1,c_2$ independent of $m,\beta$ such that 
\be
||f(u)||_{L^{\frac{s+1}{s}}}^{\frac{s+1}{s}}\leq c_1I(u)+c_2
\ee
\end{lem}
Proof:\\
If $u$ is a critical point of $I_\beta$ then $I_\beta^\prime(u)\phi=0\,\, \forall \phi\in E^m\oplus N^m$
hence
\begin{eqnarray}I_\beta(u)
& = &     I_\beta(u)-I_\beta^\prime(u)u\nonumber\\
& = &     \int_Q\frac{1}{2}uf(u)-F(u)dxdt\geq a_1(s)\int_Q|u|^{s+1}dxdt-a_2(s)
\end{eqnarray}
such constant $a_1(s),a_2(s)$ exist because $f$ satisfies (\ref{conditions}). Then we have
\begin{eqnarray}I_\beta(u)
& \geq & c_1\int_Q|f(u)|^{\frac{s+1}{s}}dxdt-c_2(s).
\end{eqnarray}
Let $u^m=w^m+v^m$ the approximate solution on $E^m\oplus N^m$ then
\be
\widehat{\Box w^m}(j,k)=\widehat{f(u^m)}(j,k)
\ee
 $\forall 2j\neq k \in E^m$, hence by lemma \ref{Riesz} and the Hausdorff-Young we have
\be
||w^m||_{C^\gamma}\leq c
\ee
with $c$ independent of $m,\beta$. Hence we can conclude that $w=\lim_{m\rightarrow+\infty}w^m\in C^\gamma$ for any $\gamma<1-\frac{s+1}{s}$.
In the following lemma we follow closely the method of \cite{R78} to get an a priori estimate on $||v||_{C^0}$ independently of $\beta$.
\begin{lem} There is a constant $c$ independent of $\beta$ such that
\be
||v(\beta)||_{C^0}\leq c
\ee
\end{lem}
Proof:\\
Let $\phi\in N$ then we have
\be
\int_Q[-\beta v_{tt}+\beta v+f(v+w)]\phi dxdt=0
\ee
or
\be
\int_Q\beta v\phi+\beta v_{t}\phi_t+[(f(v+w)-f(w))\phi]dxdt=-\int_Qf(w)\phi dxdt
\ee
and $q$ is the function defined as 
\begin{displaymath}
q(s)=
\left\{ \begin{array}{ll}
s+M            & s\geq M \\
0           & -M\leq s\leq M \\
s-M & s<M \\
\end{array}\right.
\end{displaymath}
and choose 
\be
\phi(x,t)=q(v^+(x,t))+q(v^-(x,t))
\ee
\begin{eqnarray}\int_Qv_t\phi_tdxdt
& = & \int_Qq^\prime(v^+)(v^{+}_t)^2+q^\prime(v^-)(v^{-}_t)^2+\frac{\partial}{\partial_t}(q(v^+))v^-_t+\frac{\partial}{\partial_t}(q(v^-))v^+_tdxdt\nonumber\\
& = & \int_Qq^\prime(v^+)(v^{+}_t)^2+q^\prime(v^-)(v^{-}_t)^2dxdt
\end{eqnarray}
while
\[
\int_Qv\phi dxdt\geq 0,
\]
we define
\begin{displaymath}
\psi(z)=
\left\{ \begin{array}{ll}
\min_{|\xi|\leq M_5}f(z+\xi)-f(\xi)       & z\geq 0 \\
\max_{|\xi|\leq M_5}f(z+\xi)-f(\xi)       & z<0     \\
\end{array}\right.
\end{displaymath}
which is monotone in $z$ with $\psi(0)=0$.
$Q_\delta=\{(x,t)\in Q,\,  |v(x,t)|\geq \delta \}$, $Q_{\delta}^+=\{(x,t)\in Q \,\, v(x,t)\geq \delta \}$, $Q_\delta^-=Q_\delta\setminus Q_\delta^+$.
If $v\geq 0$ then 
\begin{eqnarray}\int_{Q_\delta^+}[f(v+w)-f(w)][q^++q^-]dxdt
& \geq & \frac{\psi(\delta)}{||v|_{C^0}}\int_{{Q_\delta}^+}v(q^++q^-)dxdt\nonumber
\end{eqnarray}
Now for $v\leq -\delta$
when $v<0$ then $f(v+w)-f(w)\leq \psi(v)$ and also $q^++q^-\leq 0$
and similarly
\begin{eqnarray}\int_{Q_\delta^-}(f(v+w)-f(w)(q^++q^-)dxdt
& \geq & \frac{-\psi(-\delta)}{||v||_{C^0}}\int_{{Q_\delta}^-}v(q^++q^-)dxdt\nonumber
\end{eqnarray}
now define $\nu(z)=\min (\psi(z),\psi(-z))$ for $z\geq 0$, and $||v^\pm||_{C^0}=\max(||v^+||_{C^0},||v^-||_{C^0})$
then
\begin{eqnarray}\int_{Q_\delta}[f(v+w)-f(w)][q^++q^-]dxdt
& \geq & \frac{\nu(\delta)}{||v||_{C^0}}[\int_{Q}v^+q^++v^-q^-dxdt-\delta\int_Q|q^+|+|q^-|dxdt]\nonumber
\end{eqnarray}
and since $sq(s)\geq M|q(s)|$, then
\be
||f(w||_{C^0}\int_Q(|q^+|+|q^-|)dxdt\geq \frac{(M-\delta)\nu(\delta)}{||v||_{C^0}}\int_Q(|q^+|+|q^-|)dxdt
\ee
then for arbitrary $M<|||v^\pm||_{C^0}$ and choosing $\delta=\frac{|||v^\pm||_{C^0}}{2}$ we deduce
\be
\nu(\frac{1}{2}||v^\pm||_{C^0})\leq 4||f(w||_{C^0}
\ee
hence $||v||_{C^0}$ is bounded independently of $\beta$. 
\section{Regularity of the solution}
Here we prove that if $f\in C^{2,1}$ then the weak solution $u$ is $C^2$. Since $||v||_{C^0},||w||_{C^0}$ are bounded independetly of $\beta$, $f(u)\in C^0$. We also have
\be
(-4j^2+k^2)\widehat{w}(j,k)=\widehat{f(x,v+w)}(j,k)\,\,\ 2j\neq\pm k.
\label{perp}
\ee
Then by lemma \ref{H1} we have $w\in H^1$. Since $f$ is smooth then too $f(w+v)\in H^1$. Then (\ref{perp}) implies $w\in H^2$ and iterating once again leads to $w\in H^3$. Now going back to the original equation:
\be
-\beta v_{tt}=\Box w-f(x,u)-\beta v
\ee
and recalling that $v\in C^2$ we deduce  $v\in H^3$ which with (\ref{perp}) implies $w\in H^4$. Iterating once more implies $v^4$ then again $w\in H^5$ and $v\in H^5$.
 Thus we can differentiating with refer to $t$ in the weak sense and we have
\be
\Box w_t-\beta v_{ttt}=-f_u^\prime(x,u)(w_t+v_t)-\beta v_t
\label{perp3}
\ee
in Fourier space. We now want to get estimates independently of $\beta$ pass to the limit and find solutions of (\ref{papier}). 
Now multiplying by $v_t(\beta)$ (this is possible
 since $v\in H^1$) and integrating we have
\be
\beta(v_t,v_t)+\beta (v_{tt},v_{tt})+(f_u^\prime(x,u)v_t,v_t)=-(f_u^\prime(x,u)w_t,v_t)
\label{perp2}
\ee
and
\be
\alpha ||v_t||_{L^2}^2<(f_u^\prime(x,u)v_t,v_t) \leq -(f_u^\prime(x,u)w_t,v_t)
\ee
Now since $f_u^\prime>\alpha>0$ and $||w||_{H^1}\leq c$ with $c$ independently of $\beta$ hence there is a constant $c$ independent of $\beta$
such that $||v_t||_{L^2}\leq c$. This combined with (\ref{perp}) implies $||w||_{H^2}\leq c$ where $c$ is independent of $\beta$. Differentiating (\ref{perp3}) with refer to $t$ we get 
\be
\beta v_{tt}+\Box w_{tt}-\beta v_{tttt}+f_u^\prime(x,u)v_{tt}=-f_{uu}^{\prime\prime}(u)v_t^2-f_{uu}^{\prime\prime}(u)w_t^2-2f_{uu}^{\prime\prime}(u)w_tv_t-f_u^\prime(x,u)w_{tt}
\label{perp6}
\ee
Now we multiply (\ref{perp6}) by $v_{tt}$ and estimate the $L^2$ norm of the first term of the RHS. 
\begin{eqnarray} (f_{uu}^{\prime\prime}(x,u)v_t^2,v_{tt})
& \leq & c(f)\int_0^\pi\int_0^{2\pi}v_t^2|v_{tt}|dxdt\nonumber\\
& \leq & c(f)(\int_0^\pi\int_0^{2\pi}v_t^4dxdt)^{\frac{1}{2}}(\int_0^\pi\int_0^{2\pi}|v_{tt}|^2dxdt)\nonumber
\end{eqnarray} 
we then deduce
\begin{eqnarray} (f_{uu}^{\prime\prime}(x,u)v_t^2,v_{tt})
& \leq & c(f)||v_t||_{L^2}^{\frac{3}{2}}||v_t||_{H^1}^{\frac{3}{2}}\\
\label{GagNC}
& \leq & c(f)||v_t||_{H^1}^{\frac{3}{2}}\nonumber
\label{perp7}
\end{eqnarray}
where the constant $c(f)$ is independent of $\beta$ and the inequalities in the previous argument stems from the Gagliardo-Nirenberg inequality.
The $L^2$ norms of the terms in the RHS can be estimated by noting that $f(u)\in H^1$,$w\in C^{1,\gamma}$,$0<\gamma<\frac{1}{2}$, and that the respective norms can be estimated are independently of $\beta$:
\be
(f_{uu}^{\prime\prime}(u)w_t^2,v_{tt}) \leq c||v_{tt}||_{L^2}
\label{perp8}
\ee
\be
(2f_{uu}^{\prime\prime}(x,u)w_t,v_{tt}) \leq c||v_{tt}||_{L^2}
\label{perp9}
\ee
\be
-(f_u^\prime(x,u)w_{tt},v_{tt})\leq c||w_{tt}||_{L^2}||v_{tt}||_{L^2}
\label{perp10}
\ee
recalling (\ref{perp6}), multiplying by $v_{tt}$
\be
\beta(v_{tt},v_{tt})+\beta(v_{ttt},v_{ttt})+(f_u^\prime(x,u)v_{tt},v_{tt})=(-f_{uu}^{\prime\prime}(x,u)v_t,v_{tt})+(-f_{uu}^{\prime\prime}(x,u)w_t,v_{tt})
+(-f_u^\prime(x,u)w_{tt},v_{tt})
\label{perp5}
\ee
We can now continue from (\ref{perp6}),(\ref{perp7}),(\ref{perp8}),(\ref{perp9}),(\ref{perp10}) and we have
\begin{eqnarray}\beta(v_{tt},v_{tt})+\beta(v_{ttt},v_{ttt})+(f_u^\prime(x,u)v_{tt},v_{tt})
& \leq & c||v_t||_{H^1}^{\frac{3}{2}}
\end{eqnarray}
thus there exists $c$ independent of $\beta$ such that $||v_{tt}||_{L^2}\leq c$ where $c$ is independent of $\beta$. At this stage we can conclude that there is a constant $c$ independent of $\beta$ such that $||f(u)||_{H^2}\leq c$. Combining this with (\ref{perp}) we have $||w||_{H^3}\leq c$ with $c$ is independent of $\beta$, $w\in C^{2,\frac{1}{2}}$ and $v\in C^1$ with upper bounds independent of $\beta$. We have now proved that if $f$ is $C^2$ then the solution is $u\in H^2\cap C^1$ is a weak solution of the equation.
We now differentiate (\ref{perp6}) we have
\begin{eqnarray}\beta v_{ttt}+\Box w_{ttt}+f^\prime(u)v_{ttt}
& = & -f_{uu}^{\prime\prime}(x,u)v_{tt}-f_{uuu}^{\prime\prime\prime}(x,u)(v_t+w_t)v_t^2-f_{uu}^{\prime\prime}(x,u)2v_tw_{tt}-f_{uuu}^{\prime\prime\prime}(x,u)w_{t}^2\nonumber\\
&   & -f_{uu}^{\prime\prime}(x,u)2w_tw_{tt}-2f_{uuu}^{\prime\prime\prime}(x,u)(v_t+w_{tt})w_{t}v_t-2f_{uu}^{\prime\prime}(x,u)w_{tt}v_t\nonumber\\
&   & -2f_{uu}^{\prime\prime}(x,u)w_tv_{tt}-f_{uu}^{\prime\prime}(x,u)(v_t+w_t)w_{tt}-f_{uu}^{\prime\prime}(x,u)w_{ttt}\nonumber
\end{eqnarray}
and multiplying both sides of the preceding equality by $v_{ttt}$ and integrating we conclude that $||v_{ttt}||_{L^2}\leq c$ where $c$ is independent of $\beta$ thus $v$ is $C^2$. Now recalling the Holder regularity bootstrap and (\ref{perp}) we get $w\in C^{3,\gamma}$, $0<\gamma<\frac{1}{2}$.
\bibliographystyle{plain}

\end{document}